\documentclass[preprint,12pt,3p,times]{elsarticle}

\usepackage{enumerate,bbm}%amsthm
\usepackage[utf8]{inputenc}
\usepackage[T1]{fontenc}
\usepackage{amsmath}
\usepackage{amsfonts}
\usepackage{mathptmx}
\usepackage{amscd,amssymb}
\usepackage{graphicx}
\usepackage{microtype}
\newtheorem{thm}{Theorem}

\newtheorem{stw}[thm]{Proposition}
\newtheorem{wn}[thm]{Corollary}

\newtheorem{uw}{Remark}
\newtheorem{conj}{Conjecture}
\journal{Bulletin des Sciences Mathématiques}

\begin{document}

\begin{frontmatter}

\title{Milnor numbers in deformations of homogeneous singularities\tnoteref{grant}}

\author[afi]{Szymon Brzostowski\corref{cor1}}%\fnref{hmm}
\ead{brzosts@math.uni.lodz.pl}

\author[afi]{Tadeusz Krasiński}
\ead{krasinsk@uni.lodz.pl}

\author[afi]{Justyna Walewska}
\ead{walewska@math.uni.lodz.pl}

\tnotetext[grant]{
Szymon Brzostowski and Tadeusz Krasiński were supported by the Polish National Science Centre (NCN), Grant No. 2012/07/B/ST1/03293. Szymon Brzostowski and Justyna Walewska were supported by the Polish National Science Centre (NCN) Grant No. 2013/09/D/ST1/03701.}
\cortext[cor1]{Corresponding author}
\address[afi]{ Faculty of Mathematics and Computer Science, University of Łódź, ul. Banacha 22, 90--238 Łódź, Poland}

\begin{abstract}
Let $f_{0}$ be a plane curve singularity. We study the Minor numbers of
singularities in deformations of $f_{0}$. We completely describe the set of
these Milnor numbers for homogeneous singularities $f_{0}$ in the case of
non-degenerate deformations and obtain some partial results on this set in the
general case.
\end{abstract}


\end{frontmatter}

\renewcommand{\theenumi}{\arabic{enumi}.} \renewcommand{\labelenumi}{\theenumi}

\section{Introduction}

Let $f_{0}:(\mathbb{C}^{n},0)\rightarrow(\mathbb{C},0)$ be an isolated
singularity (in the sequel a singularity means an isolated singularity) and
$\mu(f_{0})$ its Milnor number at $0$. Consider an arbitrary holomorphic
deformation $(f_{s})_{s\in S}$ of $f_{0}$, where $s$ is a single parameter
defined in a neighborhood $S$ of $0\in\mathbb{C}$. By the semi-continuity (in
the Zariski topology) of Milnor numbers in families of singularities,
$\mu(f_{s})$ is constant for sufficiently small $s\not =0$ and $\mu(f_{s}%
)\leq\mu(f_{0})$. Denote this constant value by $\mu^{{\operatorname*{gen}}}(f_{s})$ and call it \textit{generic Milnor
number of the deformation} $(f_{s})$. Let
\[
\mathcal{M}(f_{0})=(\mu_{0}(f_{0}),\mu_{1}(f_{0}),\ldots,\mu_{k}(f_{0}))
\]
be the strictly decreasing sequence of generic Milnor numbers of all possible
deformations of $f_{0}$. In particular $\mu_{0}(f_{0})=\mu(f_{0})>\mu
_{1}(f_{0})>\ldots>\mu_{k}(f_{0})=0$. If $f_{0}$ is a fixed singularity then
the sequence $\mathcal{M}(f_{0})$ will be denoted shortly $(\mu_{0},\mu
_{1},\ldots,\mu_{k}).$ Analogously we define
\[
\mathcal{M}^{{\operatorname*{nd}}}(f_{0})=(\mu_{0}%
^{{\operatorname*{nd}}}(f_{0}),\mu_{1}^{{\operatorname*{nd}}%
}(f_{0}),\dots,\mu_{l}^{{\operatorname*{nd}}}(f_{0})),
\]
the strictly decreasing sequence of generic Milnor numbers of all possible
non-degenerate deformations of $f_{0}$ (it means that any element of the
family $(f_{s})$ is a Kouchnirenko non-degenerate singularity). Notice that
the sequence $\mathcal{M}^{{\operatorname*{nd}}}(f_{0})$ is a
subsequence of $\mathcal{M}(f_{0})$. The problem of description of
$\mathcal{M}(f_{0})$ and $\mathcal{M}^{{\operatorname*{nd}}}(f_{0})$
was posed by A. Bodin \cite{Bo} who, in turn, generalized related problems
posed by A'Campo (unpublished) and V.I. Arnold \cite{Arnold} (Problems
1975-15, 1982-12). It is a non-trivial problem because by Gusein-Zade
\cite{G-S}, see also Brzostowski-Krasi\'{n}ski \cite{BK}, Bodin \cite{Bo},
Walewska \cite{W1}, \cite{W2}, the sequence $\mathcal{M}(f_{0})$ and
consequently $\mathcal{M}^{{\operatorname*{nd}}}(f_{0})$ is not always
equal to the sequence of all non-negative integers less than $\mu(f_{0})$. For
instance, $\mathcal{M}(x^{4}+y^{4})=(9,7,6,\ldots,1,0)$, $\mathcal{M}%
^{{\operatorname*{nd}}}(x^{4}+y^{4})=(9,6,5,\ldots,1,0)$. In the paper
we will consider the class of homogeneous singularities in the plane. We
describe completely the sequence $\mathcal{M}^{{\operatorname*{nd}}%
}(f_{0})$ for homogeneous plane curve singularities $f_{0}$ and we give some
partial results on $\mathcal{M}(f_{0})$. The main results are:

\begin{thm}
\label{t:1} If $f_{0}:(\mathbb{C}^{2},0)\rightarrow(\mathbb{C},0)$ is a
homogeneous singularity of degree $d\geq2$ (it means $f_{0}$ is a homogeneous
polynomial of degree $d$ without multiple factors) then%

\begin{equation}
\mathcal{M}^{{\operatorname*{nd}}}(f_{0})=((d-1)^{2},(d-1)(d-2),\ldots
,1,0),\label{r:odd}%
\end{equation}
if $d$ is odd or $d\leq4$ or $f_{0}$ is non-convenient, and
\begin{equation}
\mathcal{M}^{{\operatorname*{nd}}}(f_{0})=((d-1)^{2},(d-1)(d-2),\ldots
,\widehat{d^{2}-4d+2},\ldots,1,0),\label{r:even}%
\end{equation}
if $d$ is even $\geq6$ and $f_{0}$ is convenient ($\widehat{a}$ means the
symbol $a$ is omitted).
\end{thm}

\begin{uw}
The value of the first jump in the above sequences $(d-1)^{2}-(d-1)(d-2)=d-1$
has been given by A. Bodin \cite{Bo}. The value of the second one equal to $1
$ was established by J. Walewska \cite{W1}, \cite{W2}.
\end{uw}

\begin{uw}
Since the Milnor number of a non-degenerate singularity depends only on its
Newton diagram (see the Kouchnirenko Theorem in Preliminaries) we obtain that
Theorem \ref{t:1} holds also for semi-homogeneous singularities i.e. for
singularities of the form $\widetilde{f_{0}}=f_{0}+g,$ where $f_{0}$ is a
homogeneous singularity of degree $d$ and $g$ is a holomorphic function of
order $>d.$
\end{uw}

In case $\mathcal{M}(f_{0})$ we can only complement the sequence
$\mathcal{M}^{{\operatorname*{nd}}}(f_{0})$ by some numbers.

\begin{thm}
\label{t:??} If $f_{0}:(\mathbb{C}^{2},0)\rightarrow(\mathbb{C},0)$ is a
homogeneous singularity of degree $d\geq3$ then
\[
\mathcal{M}(f_{0})=((d-1)^{2},\mu_{1},\ldots,\mu_{r}%
,(d-1)(d-2)+1,(d-1)(d-2),\ldots,1,0),
\]
where $\mu_{1},\ldots,\mu_{r}$ is an unknown subsequence (may be empty).
\end{thm}

For a particular homogeneous singularity we have a more precise result.

\begin{thm}
\label{t:2} For the singularity $f_{0}(x,y)=x^{d}+y^{d}$, $d\geq2,$
\begin{equation}
\mu_{1}(f_{0})=(d-1)^{2}-\left[  \dfrac{d}{2}\right]  ,\label{22}%
\end{equation}
where $[a]$ means the integer part of a real number $a$.
\end{thm}

The value of $\mu_{1}(f_{0})$ in Theorem \ref{t:2}, found for very specific
singularities, could not be extend to the whole class of homogeneous
singularities of degree $d$, because $\mathcal{M}(f_{0})$ depends on
coefficients of $f_{0}$. Precisely we have

\begin{thm}
\label{t:mu1<(d-1)^2-[d/2]} For homogeneous singularities $f_{0}$ of degree
$d\geq2$
\[
\mu_{1}(f_{0})\leq(d-1)^{2}-\left[  \dfrac{d}{2}\right]
\]
and for $f_{0}$ of degree $d\geq5$ with generic coefficients we have
\begin{equation}
\mu_{1}(f_{0})<(d-1)^{2}-\left[  \dfrac{d}{2}\right]  .\label{23}%
\end{equation}

\end{thm}

\begin{uw}
Having generic coefficients means: there is a proper algebraic subset $V$ in
the space $\mathbb{C}^{d+1}$ of coefficients of homogeneous polynomials of
degree $d$ such that for homogeneous singularities $f_{0}$ of degree $d$ with
coefficients outside $V$ the inequality \textup{(\ref{23})} holds.
\end{uw}

\section{Preliminaries}

Let $f_{0}:(\mathbb{C}^{n},0)\rightarrow(\mathbb{C},0)$ be \textit{an isolated
singularity}, i.e. $f_{0}$ is a germ of a holomorphic function having an
isolated critical point at $0\in\mathbb{C}^{n}$ and $0\in\mathbb{C}$ as the
corresponding critical value. \textit{A deformation} of $f_{0}$ is the germ of
a holomorphic function $f(z,s):(\mathbb{C}^{n}\times\mathbb{C},0)\rightarrow
(\mathbb{C},0)$ such that:

\begin{enumerate}
\item $f(z,0)=f_{0}(z)$, \label{c:2}

\item $f(0,s)=0$.
\end{enumerate}

The deformation $f(z,s)$ of the singularity $f_{0}$ will be treated as a
family $(f_{s})$ of germs, taking $f_{s}(z)=f(z,s)$. Since $f_{0}$ is an
isolated singularity, $f_{s}$ has also isolated singularities near the origin,
for sufficiently small $s$ (\cite{GLS}, Ch.I, Thm $2.6$). Then the Milnor
numbers $\mu(f_{s})$ of $f_{s}$ at $0$ are defined. Since the Milnor number is
upper semi-continuous in the Zariski topology in families of singularities
(\cite{GLS}, Prop. 2.57) there exists an open neighborhood $S$ of
$0\in\mathbb{C}$ such that

\begin{enumerate}
\item $\mu(f_{s})=\operatorname*{const}$ for $s\in S\setminus\{0\}$,

\item $\mu(f_{0})\geq\mu(f_{s})$ for $s\in S$.
\end{enumerate}

Consequently, the notions of $\mu^{{\operatorname*{gen}}}(f_{s})$,
$\mathcal{M}(f_{0})$, and $\mathcal{M}^{{\operatorname*{nd}}}(f_{0})$
in the Introduction are well-defined.

Let $\mathbb{N}$ be the set of nonnegative integers and $\mathbb{R}_{+}$ be the set of
nonnegative real numbers. Let $f_{0}(x,y)=\sum_{(i,j)\in\mathbb{N}^{2}}%
a_{ij}x^{i}y^{j}$ be a singularity. Put $\operatorname{Supp}(f_{0}%
):=\{(i,j)\in\mathbb{N}^{2}:a_{ij}\not =0\}.$ The \textit{Newton diagram} of
$f_{0}$ is defined as the convex hull of the set $\bigcup_{(i,j)\in
\operatorname{Supp}(f_{0})}\left(  (i,j)+\mathbb{R}_{+}^{2}\right)  $ and is
denoted by $\Gamma_{+}(f_{0})$. The boundary (in $\mathbb{R}^{2}$) of the
diagram $\Gamma_{+}(f_{0})$ is the sum of two half-lines and a finite number
of compact line segments. The set of those line segments will be called
\textit{the Newton polygon of} $f_{0}$ and denoted by $\Gamma(f_{0})$. For
each segments $\gamma\in\Gamma(f_{0})$ we define a weighted homogeneous
polynomial $(f_{0})_{\gamma}:=\sum_{(i,j)\in\gamma}a_{ij}x^{i}y^{j}.$ A
singularity $f_{0}$ is called \textit{non-degenerate} (in the Kouchnirenko
sense) \textit{on a segment} $\gamma\in\Gamma(f_{0})$ if and only if the
system of equations%

\[
\displaystyle{\frac{\partial(f_{0})_{\gamma}}{\partial x}}%
(x,y)=0,\;\displaystyle{\frac{\partial(f_{0})_{\gamma}}{\partial y}}(x,y)=0,
\]
has no solutions in ${\mathbb{C}}^{\ast}\times{\mathbb{C}}^{\ast}$. $f_{0}$ is
called \textit{non-degenerate} if and only if it is non-degenerate on every
segment $\gamma\in\Gamma(f_{0})$. A singularity is called \textit{convenient}
if $\Gamma_{+}(f_{0})$ intersects both coordinate axes in $\mathbb{R}^{2}$.
For such singularities we denote by $S$ the area of the domain bounded by the
coordinate axes and the Newton polygon $\Gamma(f_{0})$. Let $a$ (resp. $b$) be
the distance of the point $(0,0)$ to the intersection of $\Gamma_{+}(f_{0})$
with the horizontal (resp. vertical) axis. The number
\begin{equation}
\nu(f_{0}):=2S-a-b+1,\tag{K}\label{r:K}%
\end{equation}
is called the \textit{Newton number of the singularity} $f_{0}$. Let us recall
the Planar Kouchnirenko Theorem.

\begin{thm}[\cite{Kush}]
For a convenient singularity $f_{0}$ we have:

\begin{enumerate}
\item $\mu(f_{0})\geq\nu(f_{0})$,

\item if $f_{0}$ is non-degenerate then $\mu(f_{0})=\nu(f_{0})$.
\end{enumerate}
\end{thm}

The Newton number of singularities is monotonic with respect to the Newton
diagrams of these singularities (with the relation of inclusion).

\begin{stw}
[\cite{LenarcikManuscripta08}, \cite{Gwo08}]If $f_{0}$ and $\tilde{f_{0}}$ are
convenient singularities and $\Gamma_{+}(f_{0})\subset\Gamma_{+}(\tilde{f_{0}%
})$ then $\nu(f_{0})\geq\nu(\tilde{f_{0}})$.
\end{stw}

\begin{wn}
\label{c:monotonicznosc} If $f_{0}$ and $\tilde{f_{0}}$ are convenient,
non-degenerate singularities and $\Gamma_{+}(f_{0})\subset\Gamma_{+}%
(\tilde{f_{0}})$ then $\mu(f_{0})\geq\mu(\tilde{f_{0}})$.
\end{wn}

In the paper we will use ,,global'' results concerning projective algebraic
curves proved by A. P\l oski.

\begin{thm}[{\cite[Thm 1.1]{Ploski13}}]\label{t:ploski} Let
$f=0$, $f\in\mathbb{C}[X,Y]$, be a plane algebraic curve of degree $d>1$ with
an isolated singular point at $0\in\mathbb{C}^{2}$. Suppose that
\textup{$\operatorname*{ord}$}$_{0}(f)<d$. Then
\[
\mu(f)\leq(d-1)^{2}-\left[  \dfrac{d}{2}\right]  .
\]

\end{thm}

\begin{uw}
The assumption \textup{$\operatorname*{ord}$}$_{0}(f)<d$ in the above theorem
means that $f$ is not a homogeneous polynomial. If $f$ is a homogeneous
polynomial of degree $d$ with an isolated singular point at $0\in
\mathbb{C}^{2}$ then obviously $\mu(f)=(d-1)^{2}$.
\end{uw}

\begin{thm}[{\cite[Thm 1.4]{Ploski13}}]\label{25} Let $f$ be a
polynomial of degree $d>2$, $d\not =4$. Then the following two conditions are equivalent

\begin{enumerate}
\item The curve $f=0$ passes through the origin and $\mu_{0}(f)=(d-1)^{2}%
-\left[  d/2\right]  $,

\item The curve $f=0$ has $d-\left[  d/2\right]  $ irreducible components.
Each irreducible component of the curve passes through the origin. If
$d\equiv0$ \textup{(mod $2$)} then all components are of degree $2$ and
intersect pairwise at $0$ with multiplicity $4$. If $d\not \equiv 0$
\textup{(mod $2$)} then all but one component are of degree $2$ and intersect
pairwise at $0$ with multiplicity $4$, the remaining component is linear and
tangent to all components of degree $2$.\label{condition:2 Ploski Theorem}
\end{enumerate}
\end{thm}

\section{Proof of Theorem \ref{t:1}}

Let $f_{0}$ be a homogeneous isolated singularity of degree $d$ i.e.
\[
f_{0}(x,y)=a_{0}x^{d}+\ldots+a_{d}y^{d},a_{i}\in\mathbb{C},\;\;d\geq
2,\;\;f_{0}\not =0
\]
and $f_{0}$ has no multiple factors in $\mathbb{C}[x,y]$. Geometrically it is
an ordinary singularity of $d$ lines intersecting at the origin. Notice that
$f_{0}$ is non-degenerate.

By A. Bodin \cite{Bo} and J. Walewska \cite{W1} for any non-degenerate
deformation $(f_{s})$ of $f_{0}$, for which $\mu^{{\operatorname*{gen}%
}}(f_{s})\not =\mu(f_{0})$ we have $\mu^{{\operatorname*{gen}}}%
(f_{s})\leq\mu(f_{0})-(d-1)=(d-1)^{2}-(d-1)=(d-1)(d-2)=d^{2}-3d+2$.

\textbf{A}. Assume first that $f_{0}$ is convenient i.e. $a_{0}a_{d}\not =0$.
Since we consider only non-degenerate deformations of $f_{0}$, by the
Kouchnirenko Theorem we may assume that
\begin{equation}
f_{0}(x,y)=x^{d}+y^{d},\qquad d\geq2.\label{r:x^d+y^d}%
\end{equation}
We will apply induction with respect to the degree $d$. It is easy to find
non-degenerate deformations $(f_{s})$ of $f_{0}$ for the degrees $d=2,3,4,$
whose generic Milnor numbers realize all the numbers $\leq d^{2}-3d+2$. This
gives
\[
\mathcal{M}^{{\operatorname*{nd}}}(x^{2}+y^{2})=(1,0),
\]%
\[
\mathcal{M}^{{\operatorname*{nd}}}(x^{3}+y^{3})=(4,2,1,0),
\]%
\[
\mathcal{M}^{{\operatorname*{nd}}}(x^{4}+y^{4})=(9,6,5,\ldots,1,0),
\]
(in the last case one can use some of deformations given below). Let us
consider singularity (\ref{r:x^d+y^d}) where $d\geq5$. It is easy to check (by
the Kouchnirenko Theorem) that the deformations

\begin{enumerate}
\item $f_{s}=f_{0}+sy^{d-1}+sx^{l-d+1}y^{2d-l-2}$ for $d-1\leq l\leq2d-3$ have
generic Milnor numbers $d^{2}-4d+4,\ldots,d^{2}-3d+2,$ respectively,

\item $f_{s}=f_{0}+sy^{d-2}$ has generic Milnor number $d^{2}-4d+3$,

\item $f_{s}=f_{0}+sy^{d-2}sx^{d-[\frac{d}{2}]}y^{[\frac{d}{2}]-1}$ has
generic Milnor number $d^{2}-4d+2$ for $d$ odd and $d^{2}-4d+3$ for $d$ even.

\item $f_{s}=f_{0}+sy^{d-1}+sx^{l}y^{d-l-2}$ for $1\leq l\leq d-3$ have
generic Milnor numbers $d^{2}-5d+5$, \ldots, $d^{2}-4d+1$, respectively.
\end{enumerate}

The above deformations ,,realize'' all integers from $d^{2}-5d+5$ to
$d^{2}-3d+2$ with exception of the number $d^{2}-4d+2$ in the case $d$ is
even. Now we use induction hypothesis. Notice that for $(d-1)$ we have
$((d-1)-1)((d-1)-2)=d^{2}-5d+6>d^{2}-5d+5$. Hence, if $d$ is odd then $(d-1)$
is even and by induction hypotheses we may ,,realize'' all integers from $0$ to
$d^{2}-5d+6$ with the exception of the number $(d-1)^{2}-4(d-1)+2=d^{2}-6d+7$.
But the deformation $f_{s}=f_{0}+sy^{d-1}+sx^{d-5}y^{2}$ of $f_{0}$ has
generic Milnor number equal to $d^{2}-6d+7$. This gives formula (\ref{r:odd}).

If $d$ is even, then $(d-1)$ is odd and by induction hypothesis we may find
deformations of $f_{0}$ realizing all integers from $0$ to $d^{2}-5d+6$.
Consequently in these cases we have found deformations of $f_{0}$ realizing
all integers from $0$ to $(d-1)(d-2)$ with the exception of the number
$d^{2}-4d+2$. Now we prove that this number is not generic Milnor number of
any non-degenerate deformation of $f_{0}$. Assume to the contrary that there
exists a non-degenerate deformation $(f_{s})$ of $f_{0}(x,y)=x^{d}+y^{d}$,
$d\geq6$, $d$ even, for which
\[
\mu^{\operatorname*{gen}}(f_{s})=d^{2}-4d+2.
\]

Since for sufficiently small $s\not =0$ the Newton polygons of $f_{s}$ are the
same we consider the following cases:

\textbf{I.} ${\operatorname*{ord}}_{(x,y)}f_{s}\leq d-2$. Then there
are points $(i,j)$ in ${\operatorname{Supp}}f_{s}$, $s\not =0$, such
that $i+j\leq d-2$. Take any such point $(i,j)$. Consider subcases:

\textbf{Ia.} $(i,j)\not =(0,d-2)$ and $(i,j)\not =(d-2,0)$. Consider the
non-degenerate auxiliary deformation of $f_{0}$%

\[
\tilde{f_{s}}(x,y):=\left\{
\begin{array}
[c]{lll}%
f_{0}(x,y)+sx^{i}y^{d-2-i} & \text{ if } & i>0,\\
f_{0}(x,y)+sxy^{d-3} & \text{ if } & i=0.
\end{array}
\right.
\]
It is easy to see that $\Gamma_{+}(\tilde{f_{s}})\subset\Gamma_{+}(f_{s})$ for
$s\not =0$. Then by Corollary \ref{c:monotonicznosc}%

\[
\mu^{\operatorname*{gen}}(\tilde{f_{s}})\geq\mu^{\operatorname*{gen}}%
(f_{s})=d^{2}-4d+2.
\]
But by formula (\ref{r:K}) we obtain%

\[
\mu^{\operatorname*{gen}}(\tilde{f_{s}})=d^{2}-4d+1,
\]
a contradition.

\textbf{Ib.} $(i,j)=(0,d-2)$ or $(i,j)=(d-2,0)$. Both cases are similar, so we
will consider only the case $(i,j)=(0,d-2)$. We define the auxiliary
singularity
\[
\tilde{f_{0}}(x,y):=y^{d-2}+x^{d}.
\]

By formula (\ref{r:K}) $\mu(\tilde{f_{0}})=d^{2}-4d+3$ and obviously
$\Gamma_{+}(\tilde{f_{0}})\subset\Gamma_{+}(f_{s})$ for $s\not =0$. Hence
$\tilde{f_{s}}(x,y):=f_{s}(x,y)-y^{d}+\alpha y^{d-2}$ for some generic
$0\neq\alpha\in\mathbb{C}$, would be a non-degenerate deformation of
$\tilde{f_{0}}$ such that $\Gamma_{+}(\tilde{f_{s}})=\Gamma_{+}(f_{s})$ for
$s\not =0$. Hence $\mu^{\operatorname*{gen}}(\tilde{f_{s}})=d^{2}-4d+2$. This
gives $\mu(\tilde{f_{0}})-\mu^{\operatorname*{gen}}(\tilde{f_{s}})=1$, that is
the first jump of Milnor numbers for the singularity $\tilde{f_{0}%
}(x,y)=y^{d-2}+x^{d}$ is equal to $1$. This is impossible by Bodin result
(\cite{Bo}, Section 7) because $d\geq6$ and $\operatorname{GCD}(d-2,d)=2$ (he
proved that the first jump for non-degenerate deformations is equal to $2$ in
this case).

\textbf{II.} \textup{$\operatorname*{ord}$}$_{(x,y)}f_{s}>d-2$. Then
$\Gamma_{+}(f_{s})\subset\Gamma_{+}(x^{d-1}+y^{d-1})$. Hence $\mu
^{\operatorname*{gen}}(f_{s})\geq\mu(x^{d-1}+y^{d-1})=(d-2)^{2}=d^{2}-4d+4$,
which contradicts the supposition that $\mu^{{\operatorname*{gen}}%
}(f_{s})=d^{2}-4d+2$.

\textbf{B.} Assume now that $f_{0}$ is non-convenient i.e. $f_{0}%
(x,y)=x\tilde{f_{0}}(x,y)$ (case I) or $f_{0}(x,y)=y\tilde{f_{0}}(x,y)$ (case
II) or $f_{0}(x,y)=xy\tilde{\tilde{f_{0}}}(x,y)$ (case III), where
$\tilde{f_{0}} $ is convenient of degree $d-1$ and $\tilde{\tilde{f_{0}}}$ is
convenient of degree $d-2$. Take any integer $k\leq(d-1)(d-2)$ and consider cases:

1. $k\not =d^{2}-4d+2$ or $d$ odd. Then there exists a deformation $(f_{s}%
^{1})$ of $x^{d}+y^{d}$ such that $\mu^{{\operatorname*{gen}}}%
(f_{s}^{1})=k$. Let $f_{s}^{2}:=f_{s}^{1}-x^{d}-y^{d}$. Then for the
deformation $f_{s}:=f_{0}+f_{s}^{2}+sx^{d}$ in case I or $f_{s}:=f_{0}%
+f_{s}^{2}+sy^{d}$ in case II or $f_{s}:=f_{0}+f_{s}^{2}+sx^{d}+sy^{d}$ in
case III we obviously have $\mu^{{\operatorname*{gen}}}(f_{s})=k$.

2. $k=d^{2}-4d+2$, $d$ is even and $d\geq6$. For the deformation
$f_{s}(x,y):=f_{0}(x,y)+sy^{d-1}+sx^{2}y^{d-4}+sx^{d+2}$ in case I and III and
$f_{s}(x,y):=f_{0}(x,y)+sx^{d-1}+sy^{2}x^{d-4}+sy^{d+2}$ in case II (the
summands $sx^{d+2}$ and $sy^{d+2}$ are superfluous; they have been added in
order to use formula (\ref{r:K})) we have
\[
\mu^{{\operatorname*{gen}}}(f_{s})=d^{2}-4d+2.
\]
This ends the proof of Theorem \ref{t:1}.

\section{Proof of Theorem \ref{t:??}}

For any holomorphic function germs $f,g$ at $0\in\mathbb{C}^{2}$ by
$i_{0}(f,g)$ we will denote the \textit{intersection multiplicity} of the
plane curve singularities $f=0$ and $g=0$ at $0.$ Since $\mathcal{M}%
^{{\operatorname*{nd}}}(f_{0})$ is a subsequence of $\mathcal{M}%
(f_{0})$, it suffices to prove that the number $(d-1)(d-2)+1$ and the number
$d^{2}-4d+2$ for $d\geq6$ are generic Milnor numbers of some deformations of
$f_{0}$. Consider first the number $(d-1)(d-2)+1$. Let $f_{0}=L_{1}\cdots
L_{d}$ be a factorization of $f_{0}$ into linear forms (no pair of them are
proportional). We define the deformation of $f_{0}$ by
\[
f_{s}=sL_{1}^{d-1}+f_{0}.
\]
Take $s\not =0$. Without loss of generality we may assume that $L_{1}(x,y)=x$.
Then $f_{0}(x,y)=x(\alpha y^{d-1}+\ldots)$ where $\alpha\not =0$. Hence
\[%
\begin{split}
\mu(f_{s}) &  =\mu(sx^{d-1}+x(\alpha y^{d-1}+\ldots))=\\
&  ={i}_{0}((d-1)sx^{d-2}+(\alpha y^{d-1}+\ldots)+x\dfrac
{\partial(\alpha y^{d-1}+\ldots)}{\partial x},x((d-1)\alpha y^{d-2}%
+\ldots))=\\
&  ={i}_{0}((d-1)sx^{d-2}+(\alpha y^{d-1}+\ldots)+x\dfrac
{\partial(\alpha y^{d-1}+\ldots)}{\partial x},x)+\\
&  +{i}_{0}((d-1)sx^{d-2}+(\alpha y^{d-1}+\ldots)+x\dfrac
{\partial(\alpha y^{d-1}+\ldots)}{\partial x},(d-1)\alpha y^{d-2}+\ldots)=\\
&  =(d-1)+(d-2)(d-2)=(d-1)(d-2)+1.
\end{split}
\]

Consider now the number $d^{2}-4d+2$ for $d\geq6$. Let $f_{0}=L_{1}\cdots
L_{d}$ be a factorization of $f_{0}$ into linear forms. We may assume that
$L_{1}(x,y)=\alpha x+\beta y$ where $\alpha\not =0$. If we take a linear
change of coordinates $\Phi:x^{\prime}=L_{1}(x,y)$, $y^{\prime}=y$ in
$\mathbb{C}^{2}$ then the homogeneous singularity $\tilde{f_{0}}(x^{\prime
},y^{\prime}):=f_{0}\circ\Phi^{-1}(x^{\prime},y^{\prime})=x^{\prime}%
f_{1}(x^{\prime},y^{\prime})$ is non-convenient and degree $d$. Then by
Theorem \ref{t:1} there exists a deformation $(\tilde{f_{s}})_{s\in S}$ of
$\tilde{f_{0}}$ such that $\mu^{{\operatorname*{gen}}}(\tilde{f_{s}%
})=d^{2}-4d+2$. Hence for the deformation $f_{s}:=\tilde{f_{s}}\circ\Phi$,
$s\in S,$ of $f_{0}$ we obtain $\mu(f_{s})=\mu(\tilde{f_{s}}\circ\Phi
)=\mu(\tilde{f_{s}})=d^{2}-4d+2$ for $s\not =0$. Then $\mu
^{{\operatorname*{gen}}}(f_{s})=d^{2}-4d+2$.

\section{Proof of Theorem \ref{t:2}}

Let $f_{0}(x,y)=x^{d}+y^{d}$, $d\geq2$. Let us take a deformation $(f_{s})$ of
$f_{0}$ which realizes the generic Milnor number $\mu_{1}$ of $f_{0}$ i.e.%

\begin{equation}
\mu^{{\operatorname*{gen}}}(f_{s})<\mu(f_{0})\label{r:mu(f_s)<mu(f_0)}%
\end{equation}
and
\begin{equation}
\mu(f_{0})-\mu^{{\operatorname*{gen}}}(f_{s}%
)\label{r:mu(f_0)-mu(f_s) is minimal}%
\end{equation}
is minimal non-zero integer among all deformation of $f_{0}$. In order to
apply P\l {}oski Theorem \ref{t:ploski} to elements $f_{s}$, $s\not =0,$ of
the family $(f_{s})$ we have to fulfill the assumptions of this theorem. We
will achieve this by modifying the deformation $(f_{s})$ to another one
$(\tilde{f_{s}})$ which satisfies all the requested conditions. The first step
is to reduce holomorphic $f_{s}$ to polynomials (in variables $x,y$). Notice
$(f_{s})$, $s\not =0,$ is a $\mu$-constant family. So, if we omit in
$f_{s}(x,y)$ all the terms of order $>\mu^{{\operatorname*{gen}}}%
(f_{s})+1$ then we obtain a deformation $(\tilde{f_{s}})$ of $f_{0}$ such
that
\[
{\operatorname*{ord}}(f_{s}-\tilde{f_{s}})>\mu
^{{\operatorname*{gen}}}(f_{s})+1=\mu(f_{s})+1,\qquad s\not =0.
\]
Hence by well-known theorem (\cite{Arnold}, Prop. 1 and 2 in Section 5.5,
\cite{ploski85}, Prop. 1.2 and Lemma 1.4)
\[
\mu(f_{s})=\mu(\tilde{f_{s}})\qquad\text{for}\quad s\not =0.
\]
This implies $\mu^{{\operatorname*{gen}}}(f_{s})=\mu
^{{\operatorname*{gen}}}(\tilde{f_{s}})$. By this step we may assume in
the sequel that the deformation $(f_{s})$ of $f_{0}$ which realizes $\mu_{1}$
consists of polynomials.

The second step is to reduce the degree of $f_{s}$ to $d$. For this we apply
the method of Gabrielov and Kouchnirenko \cite{gabrielov+kush}. Notice first
that there are terms in $f_{s}$ of order $<d$ with non-zero coefficients. In
fact if \textup{$\operatorname*{ord}$}$(f_{s})\geq d$ then $\mu(f_{s})\geq
\mu(x^{d}+y^{d})=(d-1)^{2}$, a contradiction. Let
\[
f_{s}(x,y)=x^{d}+y^{d}+u_{1}(s)x^{\alpha_{1}}y^{\beta_{1}}+\ldots
+u_{k}(s)x^{\alpha_{k}}y^{\beta_{k}},
\]
where $u_{i}(s),$ $i=1,\ldots,k,$ are non-zero holomorphic functions in a
neighbourhood of $0\in\mathbb{C},$ $u_{i}(0)=0.$ Denote $d_{i}:=\alpha
_{i}+\beta_{i}$, $\gamma_{i}:=$\textup{$\operatorname*{ord}$}$(u_{i})>0$,
$u_{i}(s)=a_{i}s^{\gamma_{i}}+\ldots,$ $a_{i}\neq0.$

By the above there exist $d_{i}$ for which $d_{i}<d$. Let
$N:={\operatorname*{LCM}}(\gamma_{i})$, $v:=\max\left(  \dfrac
{N\cdot(d-d_{i})}{\gamma_{i}}\right)  $. Then $v>0$. We define a new
holomorphic deformation of $f_{0}$ depending on two parameters%
\begin{align*}
f_{s,t}(x,y) &  :=\frac{f_{t^{v}s}(t^{N}x,t^{N}y)}{t^{Nd}}\\
&  =x^{d}+y^{d}+(a_{1}s^{\gamma_{1}}+t\widetilde{u}_{1}(s,t))t^{v\gamma
_{1}+N(d_{1}-d)}x^{\alpha_{1}}y^{\beta_{1}}+\ldots
\end{align*}
for some holomorphic functions $\tilde{u_{i}}(s,t)$, $i=1,\ldots,k$. By
semi-continuity of Milnor numbers in families of singularities we obtain that
for any fixed $s\in S$, $s\not =0,$
\[
\mu(f_{s,0})\geq\mu(f_{s,t})\qquad\text{for sufficiently small}\;t.
\]
But for any fixed $s,t\neq0$ sufficiently small%

\[
\mu(f_{s,t})=\mu\left(  \dfrac{f_{t^{v}s}(t^{N}x,t^{N}y)}{t^{Nd}}\right)
=\mu(f_{t^{v}s}(x,y))=\mu_{1}.
\]

Hence $\mu(f_{s,0})\geq\mu_{1}$. But
\[
f_{s,0}(x,y)=x^{d}+y^{d}+\sum_{\substack{j\\v\gamma_{j}+N(d_{j}-d)=0}%
}a_{j}s^{\gamma_{j}}x^{\alpha_{j}}y^{\beta_{j}}.
\]
Of course $d_{j}<d$ for $j$ satisfying $v\gamma_{j}+N(d_{j}-d)=0.$ So, we have
obtained a new deformation $\tilde{f_{s}}:=f_{s,0}$ of $f_{0}$ for which
\textup{$\deg$}$\tilde{f_{s}}=d$, \textup{$\operatorname*{ord}$}$(\tilde
{f_{s}})<d$ and $\mu^{{\operatorname*{gen}}}(\tilde{f_{s}})\geq\mu_{1}%
$. By definition of $\mu_{1}$ we have either $\mu^{{\operatorname*{gen}%
}}(\tilde{f_{s}})=\mu_{1}$ or $\mu^{{\operatorname*{gen}}}(\tilde
{f_{s}})=\mu(f_{0})$. The latter case is impossible because then
$\tilde{(f_{s})}$ for $s\in S$ would be a $\mu$-constant family in which one
element is equal to $f_{0}$. Since it is a family of plane curve
singularities, the orders of this singularities are the same. Hence
\[
{\operatorname*{ord}}(f_{0})={\operatorname*{ord}}(\tilde{f_{s}%
})\qquad\text{for \ }s\in S,
\]
which is impossible.

Summing up, we have obtained a deformation $(\tilde{f_{s}})$ of $f_{0}$ for
which $\mu^{{\operatorname*{gen}}}(\tilde{f_{s}})=\mu_{1}$,
\textup{$\deg$}$(\tilde{f_{s}})=d$, \textup{$\operatorname*{ord}$}%
$(\tilde{f_{s}})<d$, $s\not =0,$ and moreover $d$-th homogeneous component of
$\tilde{f_{s}}$ is equal to $f_{0}$. So, for any fixed $s\not =0$ $f_{s}$
satisfies the assumption of Theorem \ref{t:ploski}. By this theorem
\[
\mu_{1}=\mu(\tilde{f_{s}})\leq(d-1)^{2}-\left[  \dfrac{d}{2}\right]  .
\]

Now we prove the opposite inequality
\begin{equation}
\mu_{1}\geq(d-1)^{2}-\left[  \dfrac{d}{2}\right]  .
\end{equation}

It suffices to give deformations $(f_{s})$ of $f_{0}$ for which $\mu
^{{\operatorname*{gen}}}(f_{s})=(d-1)^{2}-\left[  d/2\right]  .$
Consider two cases:

1. $d$ is even i.e. $d=2k$, $k\geq1$. We define a deformation of $f_{0}$ by%

\[
f_{s}(x,y):=x^{d}+(y^{2}+sx)^{\frac{d}{2}}=x^{2k}+(y^{2}+sx)^{k}.
\]
Then we easily find for $s\not =0$
\[
\mu(f_{s})={i}_{0}(2kx^{2k-1}+ks(y^{2}+sx)^{k-1},2ky(y^{2}%
+sx)^{k-1})=4k^{2}-5k+1=(d-1)^{2}-\left[  \dfrac{d}{2}\right]  .
\]

2. $d$ is odd i.e. $d=2k+1$, $k\geq1$. We define a deformation of $f_{0}$ by
\begin{equation}
f_{s}(x,y):=(x+y)\prod\limits_{i=1}^{k}(sw_{i}(x+y)+x^{2}+2\mathfrak{Re}%
(\varepsilon^{i})xy+y^{2}),\label{r:d is odd}%
\end{equation}
where $\varepsilon={\exp(\frac{2\pi i}{2k+1})}$ is a primitive root of
unity of degree $2k+1$ and $w_{i}=\dfrac{1-\mathfrak{Re}(\varepsilon^{i}%
)}{1-\mathfrak{Re}(\varepsilon)}$, $i=1,\ldots,k$. We easily check that
$(f_{s})$ is a deformation of $f_{0}(x,y)=x^{2k+1}+y^{2k+1}$. Moreover if we
denote by $\tilde{f_{0}}$, $\tilde{f_{1}}$, \ldots,$\tilde{f_{k}}$ the
successive factors in (\ref{r:d is odd}) we easily compute that ${i}%
_{0}(\tilde{f_{0}},\tilde{f_{i}})=2$, $i=1,\ldots,k,$ and ${i}%
_{0}(\tilde{f_{i}},\tilde{f_{j}})=4$, $i,j=1,\ldots,k$, $i\not =j$. Hence by a
well-known formula for the Milnor number of a product of singularities we
obtain
\[
\mu(\tilde{f_{0}}\cdots\tilde{f_{k}})=\sum\limits_{i=0}^{k}\mu(\tilde{f_{i}%
})+2\sum_{\substack{i,j=0\\i<j}}^{k}{i}_{0}(\tilde{f_{i}},\tilde
{f_{j}})-k=4k^{2}-k=(d-1)^{2}-\left[  \dfrac{d}{2}\right]  .
\]

\section{Proof of Theorem \ref{t:mu1<(d-1)^2-[d/2]}}

Let us begin with a remark. The second part of Theorem
\ref{t:mu1<(d-1)^2-[d/2]} concerns only homogeneous singularities of degree
$d\geq5.$ For degrees $d=2,3,4$ the sequences $\mathcal{M}(f_{0})$ do not
depend on the coefficients of $f_{0}$ and they are as follows%
\[
\mathcal{M}(f_{0})=\left\{
\begin{array}
[c]{lll}%
(1,0) & \text{\textup{for}} & d=2,\\
(4,3,2,1,0) & \text{\textup{for}} & d=3,\\
(9,7,6,5,4,3,2,1,0) & \text{\textup{for}} & d=4.
\end{array}
\right.
\]
For $d=2,3$ it is an easy fact and for $d=4$ it follows from \cite{BK}.
Moreover, by Theorem \ref{t:mu1<(d-1)^2-[d/2]} we obtain that for $d=5$ we
have only two possibilities
\[
\mathcal{M}(f_{0})=\left\{
\begin{array}
[c]{ll}%
(16,13,12,\ldots,1,0) & \text{\textup{for }}f_{0}\text{ with \textup{generic
coefficients}},\\
(16,14,13,12,\ldots,1,0) & \text{\textup{otherwise}.}%
\end{array}
\right.
\]

Now we may pass to the proof of Theorem \ref{t:mu1<(d-1)^2-[d/2]}. For the
first part of the theorem we repeat the reasoning in the proof of inequality
$\mu_{1}(x^{d}+y^{d})\leq(d-1)^{2}-\left[  d/2\right]  $ in Theorem \ref{t:2}
because it works for any homogeneous singularity.

For the second part of the theorem let $f_{0}(x,y)=c_{0}y^{d}+c_{1}%
y^{d-1}x+\cdots+c_{d}x^{d}$ be an arbitrary homogeneous singularity of degree
$d\geq5$. We may assume that $c_{0}=1$ (because singularities for which
$c_{0}=0$ are ,,not generic''). Denote
\[
f_{0}^{\boldsymbol{c}}(x,y):=y^{d}+c_{1}y^{d-1}x+\cdots+c_{d}x^{d}%
,\qquad\boldsymbol{c}=(c_{1},\ldots,c_{d})\in\mathbb{C}^{d}%
.\]
By the first part of the theorem $\mu_{1}(f_{0}^{\boldsymbol{c}})\leq
(d-1)^{2}-\left[  d/2\right]  $. We will find a polynomial $F(c_{1}%
,\ldots,c_{d})$, $F\not =0$, such that if
\begin{equation}
\mu_{1}(f_{0}^{\boldsymbol{c}})=(d-1)^{2}-\left[  \dfrac{d}{2}\right]
\label{r:mu1(f_0^c)=(d-1)^2-[d/2]}%
\end{equation}
then $F(\boldsymbol{c})=0$. This will give the second part of the theorem and
finish the proof.

Let us take an arbitrary $f_{0}^{\boldsymbol{c}}$ for which
(\ref{r:mu1(f_0^c)=(d-1)^2-[d/2]}) holds. Let
\[
f_{0}^{\boldsymbol{c}}(x,y)=(y-a_{1}x)\cdots(y-a_{d}x),\quad a_{i}\not =%
a_{j}\;\text{for}\;i\not =j
\]
be the factorization of $f_{0}^{\boldsymbol{c}}$ into linear parts. We will
also denote this polynomial by $f_{0}^{\boldsymbol{a}}$, where $\boldsymbol{a}%
=(a_{1},\ldots,a_{d})$ and call $a_{1},\ldots,a_{d}$ the roots of
$f_{0}^{\boldsymbol{a}}$. By the Vieta formulas connecting $c_{i}$ with
$a_{j}$ it suffices to find a polynomial $G(a_{1},\ldots,a_{d})$, $G\not =0$,
such that if $\mu_{1}(f_{0}^{\boldsymbol{a}})=(d-1)^{2}-\left[  d/2\right]  $,
then $G(\boldsymbol{a})=0$. So, take $f_{0}^{\boldsymbol{a}}$ for which
$\mu_{1}(f_{0}^{\boldsymbol{a}})=(d-1)^{2}-\left[  d/2\right]  $. Using the
same method as in the proof of Theorem \ref{t:2} there exists a deformation
$(f_{s})$ of $f_{0}^{\boldsymbol{a}}$ such that for $s\not =0$ sufficiently small

1. \textup{$\deg$}$f_{s}=d,$

2. $d$-th homogeneous component of $f_{s}$ is equal to $f_{0}^{\boldsymbol{a}%
},$

3. $\mu^{{\operatorname*{gen}}}(f_{s})=(d-1)^{2}-\left[  d/2\right]  $.

Let us fix $s\not =0$. By Theorem \ref{25} we obtain a factorization
\[
f_{s}=LQ_{1}\cdots Q_{\left[  d/2\right]  },
\]
where $L$ is either a linear form (if $d$ is odd) or $L=1$ (if $d$ is even),
$Q_{i}$ are irreducible polynomials of degree $2$, $L$ has a common tangent
with each $Q_{i}$ (if $d$ is odd) and ${i}_{0}(Q_{i},Q_{j})=4$ for
$i\not =j$. Since $Q_{i}$ is irreducible, $Q_{i}=L_{i}+\tilde{Q}_{i}$ where
$L_{i}$ is a non-zero linear form and $\tilde{Q_{i}}$ is a non-zero quadratic
form. Moreover, the equality ${i}(Q_{i},Q_{j})=4$ for $i\not =j$ implies
the all $L_{i}$ are proportional. Additionally in the odd case $L$ is also
proportional to $L_{i}$. Notice also that $d$-th homogeneous component of
$LQ_{1}\ldots Q_{\left[  d/2\right]  }$ is equal to $L\tilde{Q}_{1}%
\ldots\tilde{Q}_{\left[  d/2\right]  }$, which by above condition 2 implies
$f_{0}^{\boldsymbol{a}}=L\tilde{Q}_{1}\ldots\tilde{Q}_{\left[  d/2\right]  }$.
Now we consider the cases:

I. $d$ is odd. By renumbering $a_{1},\ldots,a_{d}$ we may assume that
$L=(y-a_{1}x)$, $\tilde{Q}_{1}=(y-a_{2}x)(y-a_{3}x),\ldots,\tilde{Q}_{\left[
d/2\right]  }=(y-a_{d-1}x)(y-a_{d}x)$. Let $L_{i}=w_{i}L$, $w_{i}\in
\mathbb{C}\setminus\{0\}$, $i=1,\ldots,\left[  d/2\right]  $. Then from the
condition ${i}_{0}(Q_{i},Q_{j})=4$, $i\not =j,$ we obtain in particular
for $i=1$, $j=2$
\[
4={i}_{0}(w_{1}L+\tilde{Q}_{1},w_{2}L+\tilde{Q}_{2})={i}%
_{0}(\tilde{Q}_{1}-\dfrac{w_{1}}{w_{2}}\tilde{Q}_{2},w_{2}L+\tilde{Q}_{2}).
\]
Since $\tilde{Q}_{1}-\dfrac{w_{1}}{w_{2}}\tilde{Q}_{2}$ is a form of degree
$2$, we get $\tilde{Q}_{1}-\dfrac{w_{1}}{w_{2}}\tilde{Q}_{2}=uL^{2}$ for some
$u\in\mathbb{C}\setminus\{0\}$. Hence there exist $z_{1},z_{2}\in
\mathbb{C}\setminus\{0\}$ such that
\[
z_{1}\tilde{Q}_{1}+z_{2}\tilde{Q}_{2}=L^{2}.
\]
Then the non-zero $z_{1},z_{2}$ satisfy the system of equations
\[%
\begin{array}
[c]{l}%
z_{1}+z_{2}=1,\\
(a_{2}+a_{3})z_{1}+(a_{4}+a_{5})z_{2}=2a_{1},\\
a_{2}a_{3}z_{1}+a_{4}a_{5}z_{2}=a_{1}^{2}.
\end{array}
\]
Hence
\[
\left\vert
\begin{array}
[c]{lll}%
1 & 1 & 1\\
a_{2}+a_{3} & a_{4}+a_{5} & 2a_{1}\\
a_{2}a_{3} & a_{4}a_{5} & a_{1}^{2}%
\end{array}
\right\vert =0.
\]
We have obtained a non-trivial relation $\tilde{G}(a_{1},\ldots,a_{5})=0$
between the roots $a_{1},\ldots,a_{5}$. So if we put
\[
G(a_{1},\ldots,a_{d}):=\prod\limits_{\sigma\in V_{d}^{5}}\tilde{G}%
(a_{\sigma(1)},\ldots,a_{\sigma(5)}),
\]
where $V_{d}^{k}$ denotes the set of all partial permutations of length $k$
from a $d$-set, then $G$ is a non-zero polynomial in $\mathbb{C}[a_{1}%
,\ldots,a_{d}]$ such that if $G(\boldsymbol{a})\not =0$ then $\mu_{1}%
(f_{0}^{\boldsymbol{a}})<(d-1)^{2}-\left[  d/2\right]  $ . This ends the proof
of the theorem in this case.

II. $d$ is even. Then $d\geq6$. By renumbering $a_{0},\ldots,a_{d}$ we may
assume that $\tilde{Q}_{1}=(y-a_{1}x)(y-a_{2}x),\ldots,\tilde{Q}_{\left[
d/2\right]  }=(y-a_{d-1}x)(y-a_{d}x)$. Let $L_{i}=w_{i}L$, $w_{i}\in
\mathbb{C}\setminus\{0\}$, $i=1,\ldots,\left[  d/2\right]  $, where $L$ is a
fixed non-zero linear form. Repeating the reasoning as in I for the equality
${i}_{0}(Q_{1},Q_{2})=4$ we get that there exist $z_{1},z_{2}\not =0$
such that
\[
z_{1}\tilde{Q}_{1}+z_{2}\tilde{Q}_{2}=L^{2}.
\]
We may assume that either $L=x$ or $L=y$ or $L=y-\alpha x$, $\alpha\not =0$.
In the first two cases we easily obtain, as in I, non-trivial relations
$G_{1}(a_{1},a_{2},a_{3},a_{4})=0$ and $G_{2}(a_{1},a_{2},a_{3},a_{4})=0$,
respectively, between the roots $a_{1},\ldots,a_{4}$. In the third case we
obtain the relation
\[
((a_{3}+a_{4})-(a_{1}+a_{2}))\alpha^{2}-2(a_{3}a_{4}-a_{1}a_{2})\alpha
+(a_{1}+a_{2})a_{3}a_{4}-a_{1}a_{2}(a_{3}+a_{4})=0
\]
between $a_{1},a_{2},a_{3},a_{4}$ and $\alpha$. But if we apply the same
reasoning to the equality ${i}_{0}(Q_{2},Q_{3})=4$ we obtain a second
relation
\[
((a_{5}+a_{6})-(a_{3}+a_{4}))\alpha^{2}-2(a_{5}a_{6}-a_{3}a_{4})\alpha
+(a_{3}+a_{4})a_{5}a_{6}-a_{3}a_{4}(a_{5}+a_{6})=0
\]
between $a_{3},a_{4},a_{5},a_{6}$ and the same $\alpha$. Hence the resultant
of these two polynomials with respect to $\alpha$ must be equal to $0,$ which
gives a non-trivial relation $G_{3}(a_{1},\ldots,a_{6})=0$ between roots
$a_{1},\ldots,a_{6}$. Hence
\[
G(a_{1},\ldots,a_{d}):=\prod\limits_{\sigma\in V_{d}^{4}}\left(  G_{1}%
G_{2}\right)  (a_{\sigma(1)},\ldots,a_{\sigma(4)})\prod\limits_{\sigma\in
V_{d}^{6}}G_{3}(a_{\sigma(1)},\ldots,a_{\sigma(6)})
\]
is a non-trivial polynomial such that if $G(\boldsymbol{a})\not =0$ then
$\mu_{1}(f_{0}^{\boldsymbol{a}})<(d-1)^{2}-\left[  d/2\right]  $. This ends
the proof of the theorem in this case.

\section{Concluding remarks}

We have completely solved the problem of possible generic Milnor numbers of
all non-degenerate deformations of homogeneous plane singularities. The same
problem for the family of all deformations is more complicated. In the
particular case $f_{0}(x,y)=x^{d}+y^{d}$, $d\geq2$, we have only found
$\mu_{1}=(d-1)^{2}-\left[  d/2\right]  $. For generic homogeneous
singularities of degree $d$ this is not longer true by Theorem
\ref{t:mu1<(d-1)^2-[d/2]}. We do not know the exact value of $\mu_{1}$ in this
generic case. We only conjecture that for generic homogeneous singularities of
degree $d$ $\mu_{1}=(d-1)^{2}-(d-2)$. If it is true then by Theorem \ref{t:??}
we would get the whole sequence $\mathcal{M}(f_{0})$ in this case.

\begin{conj}
If $f_{0}$ is a homogeneous singularity of degree $d$ with generic
coefficients then
\[
\mathcal{M}(f_{0})=((d-1)^{2},(d-1)^{2}-(d-2),\ldots,1,0).
\]
\end{conj}

\section{Acknowledgements}
The authors thank very much Arkadiusz P\l oski and participants of the seminar
,,Analytic and Algebraic Geometry'' for valuable remarks which considerably
improved the first proofs of the theorems.
%\end{acknowledgements}

\bibliographystyle{elsarticle-num}
\bibliography{bibliografia}

\end{document}